\documentclass[12pt]{article}
\usepackage[intlimits]{amsmath}
\usepackage{amssymb}
\usepackage{mathrsfs}
\usepackage[dvips]{graphicx}
\usepackage{longtable}

\title{On the chromatic numbers of spheres in $ {\mathbb R}^n $\footnote{This work is done 
under the financial support of the following grants: the grant 09-01-00294 of 
Russian Foundation for Basic Research,
the grant MD-8390.2010.1 of the Russian President, the grant NSh-691.2008.1 supporting Leading scientific schools of Russia, 
a grant of Dynastia foundation.}}

\author{A.M. Raigorodskii}

\date{}

\begin{document}

\renewcommand{\baselinestretch}{1.5}
\maketitle

\begin{large}

\section{Introduction}

In this paper, we study a classical problem going back to H. Hadwiger, E. Nelson, and P. Erd\H{o}s. Let $ (X,\rho) $ be a metric space. 
Consider a set $ {\cal A} $ of distinct positive reals. We call the value 
$$
\chi((X,\rho);{\cal A}) = \min \left\{\chi:~ X = X_1 \bigsqcup \ldots \bigsqcup X_{\chi}, ~ \forall ~ i ~ \forall ~ x, y \in X_i ~~
\rho(x,y) \not \in {\cal A}\right\}
$$
{\it the chromatic number of the space $ (X,\rho) $ with the set of forbidden distances $ {\cal A} $}. In other words, 
$ \chi((X,\rho);{\cal A}) $ is the minimum number of colours needed to paint all the points in $ X $ so that any two points at a distance 
from $ {\cal A} $ apart receive different colours. 

Various metric spaces and sets of forbidden distances have been considered by many authors. Let us briefly review the most important cases.  

\begin{enumerate}

\item{$ (X,\rho) = ({\mathbb R}^n, l_2) $, $ {\cal A} = \{1\} $. Here 
$$ 
l_2({\bf x},{\bf y}) = \sqrt{(x_1-y_1)^2 + \ldots + (x_n-y_n)^2},
$$
$$
{\bf x} = (x_1, \dots, x_n), ~~~~ {\bf y} = (y_1, \dots, y_n).
$$
This is the classical case, which is deeply investigated. We use a simpler standard notation $ \chi({\mathbb R}^n) $ for the corresponding 
chromatic number. Numerous results concerning $ \chi({\mathbb R}^n) $ can be found in the books \cite{Soi}, \cite{BMP} and surveys \cite{Szek}, 
\cite{Rai1}. For our further purposes, only the following bounds will be useful:
$$
\chi({\mathbb R}^n) \ge (\zeta_1 + o(1))^n, ~~~ \zeta_1 = \frac{1+\sqrt{2}}{2} = 1.207... ~~~ ({\rm see} ~ \cite{FW}),
$$ 
$$
\chi({\mathbb R}^n) \ge (\zeta_2 + o(1))^n, ~~~ \zeta_2 = 1.239... ~~~ ({\rm see} ~ \cite{Rai2}),
$$ 
$$
\chi({\mathbb R}^n) \le (3 + o(1))^n ~~~ ({\rm see} ~ \cite{LR}).
$$} 

\item{$ (X,\rho) = ({\mathbb R}^n, l_2) $, $ |{\cal A}| = k $, $ k \in {\mathbb N} $. Here the best known results are given in the 
paper $ \cite{Rai3} $.}

\item{$ (X,\rho) = ({\mathbb R}^n, l_2) $, $ |{\cal A}| = \infty $. Here the paper \cite{Rai4} should be cited.} 

\item{$ (X,\rho) = ({\mathbb R}^n, l_p) $, $ |{\cal A}| = k $, $ k \in {\mathbb N} $, where
$$ 
l_p({\bf x},{\bf y}) = \sqrt[p]{|x_1-y_1|^p + \ldots + |x_n-y_n|^p}, ~~~ p \in [1,\infty),
$$
$$   
l_{\infty}({\bf x},{\bf y}) = \max_{i=1, \ldots, n} |x_i-y_i|.
$$
These cases were studied in \cite{KF}, \cite{Rai5}, \cite{Rai6}, \cite{Rai7}.} 

\item{$ (X,\rho) = ({\mathbb Q}^n, l_p) $, $ |{\cal A}| = k $, $ k \in {\mathbb N} $. See \cite{Rai1}, \cite{Rai7}, \cite{Wood}, \cite{BP}
for multiple references.}

\end{enumerate}

Another interesting series of metric spaces is generated by spheres $ S_r^{n-1} $ of radii $ r \ge \frac{1}{2} $ in $ {\mathbb R}^n $: 
$ (X,\rho) = (S_r^{n-1},l_2) $, $ {\cal A} = \{1\} $. Studying 
$$
\chi(S_r^{n-1}) = \chi((S_r^{n-1},l_2);\{1\}) 
$$
was proposed by Erd\H{o}s who conjectured in \cite{E} that $ \chi(S_r^{n-1}) \to \infty $ for any fixed value of $ r > \frac{1}{2} $. 
It is obvious that $ \chi(S_{1/2}^{n-1}) = 2 $, and L. Lov\'asz proved Erd\H{o}s' conjecture in \cite{Lov} using topological tools
(see also \cite{Mat}). The exact assertion of Lov\'asz is as follows: {\it for any $ r > \frac{1}{2} $ and $ n \in {\mathbb N} $, 
the inequality holds $ \chi(S_r^{n-1}) \ge n $; if $ r < \sqrt{\frac{n}{2n+2}} \sim \frac{1}{\sqrt{2}} $, i.e., the length of any side 
of a regular $ n $-simplex inscribed into $ S_r^{n-1} $ is smaller than 1, then $ \chi(S_r^{n-1}) \le n+1 $.} Although this result is 
widely cited (see, e.g., \cite{Szek}), its second part is completely wrong (see Section 5). Actually, for every $ r > \frac{1}{2} $, the quantity 
$ \chi(S_r^{n-1}) $ grows exponentially, not linearly. 

In this paper, we will do a careful analysis of the asymptotic behaviour of the value $ \chi(S_r^{n-1}) $. We will study even some cases when 
$ r $ may depend on $ n $. 

Before proceeding to formulating our main reults, let us mention some more papers on the chromatic numbers of spheres: \cite{Sim1}, \cite{Sim2}.

\section{Statements of the main results}

The starting point for our investigation is the following assertion. 

\vskip+0.2cm

\noindent {\bf Theorem 1.} {\it For any $ r > \frac{1}{2} $, there exist a constant $ \gamma = \gamma(r) > 1 $ and a function
$ \varphi(n) = \varphi(n,r) = o(1) $, $ n \to \infty $, such that for every $ n \in {\mathbb N} $, the inequality holds 
$$
\chi(S_r^{n-1}) \ge (\gamma + \varphi(n))^n.
$$}  

\vskip+0.2cm

Theorem 1 sais that, for any fixed radius, the chromatic number grows exponentially in the dimension. Of course it is possible 
to make the value of $ \gamma(r) $ a bit more concrete. The first step in this direction is given in Theorem 2.   

\vskip+0.2cm

\noindent {\bf Theorem 2.} {\it For any $ r \in \left(\frac{1}{2}, \frac{1}{\sqrt{2}}\right) $, 
there exists a function $ \delta(n) = \delta(n,r) = o(1) $, $ n \to \infty $, such that for every
$ n \in {\mathbb N} $, the inequality holds 
$$ 
\chi(S_r^{n-1}) \ge \left(2 \left(\frac{1}{8r^2}\right)^{\frac{1}{8r^2}} 
\left(1-\frac{1}{8r^2}\right)^{1-\frac{1}{8r^2}}+\delta(n)\right)^n.
$$} 

\vskip+0.2cm

Looking at Theorem 2, we see that if $ r $ becomes closer and closer to $ \frac{1}{\sqrt{2}} $, then the constant 
$$
\gamma = 2 \left(\frac{1}{8r^2}\right)^{\frac{1}{8r^2}} 
\left(1-\frac{1}{8r^2}\right)^{1-\frac{1}{8r^2}}
$$
approaches the value $ \zeta_3 = 1.139... $ Since $ S_r^{n-1} \subset {\mathbb R}^n $ leading to 
$ \chi(S_r^{n-1}) \le \chi({\mathbb R}^n) $, one may not expect that $ \zeta_3 $ could be 
somehow replaced by anything greater than $ \zeta_2 $ (cf. Introduction). However, there is some room to spare here, and in Section 7
we will exhibit a further optimization process providing even larger constants. 

At the same time, if $ r \ge \frac{1}{\sqrt{2}} $, then we certainly have
$$
\chi(S_r^n) \ge \chi(S_{r'}^{n-1}) \ge (1.139 + o(1))^n, ~~~ r' < \frac{1}{\sqrt{2}} \le r.
$$

So, once again, 
for {\it any} fixed value of radius, the chromatic number is essentially exponential in $ n $. Comparing our results with those due to 
Lov\'asz, we get the following assertion. 

\vskip+0.2cm

\noindent {\bf Theorem 3.} {\it For any $ r > \frac{1}{2} $,  
there exists an $ n_0 $ such that for every $ n \ge n_0 $, $ \chi(S_r^{n-1}) > n+1 $.}

\vskip+0.2cm

On the one hand, Theorem 3 shows that the bound $ \chi(S_r^{n-1}) \le n+1 $ is false, provided we fix $ r $ and let $ n $ go to infinity. 
On the other hand, the result of Theorem 3 is much stronger than that of Lov\'asz only for the values of $ n $ which are big enough. So 
in small dimensions, the {\it lower} estimate $ \chi(S_r^{n-1}) \ge n+1 $ is still the best known (and true). 

The gap between exponents and linear functions is quite large. Thus, one may expect that superlinear lower bounds for $ \chi(S_r^{n-1}) $
would be possible not only for a constant $ r > \frac{1}{2} $, but also for some sequences $ r_n \to \frac{1}{2} $. The most general 
assertion of this kind is in Theorem 4. 

\vskip+0.2cm

\noindent {\bf Theorem 4.} {\it Let $ {\mathbb P} $ be the set of prime numbers. Let $ f(x) $ be such a function that for any $ x \in {\mathbb R} $,
$ x \ge 0 $, 
$$
x+f(x)=\min\{p \in {\mathbb P}:~ p > x\}.
$$
Let 
$$
m(x) = \max\{m < x:~ m \equiv 0 \pmod 4\}.
$$
Consider a sequence $ \{r_n\}_{n=1}^{\infty} $, where $ r_n > \frac{1}{2} $ for each $ n \in {\mathbb N} $. Set 
$$
p(n) = \frac{m(n)}{8r_n^2} + f\left(\frac{m(n)}{8r_n^2}\right).
$$
If
$$
\frac{m(n)}{4} < p(n) \le \frac{m(n)}{2}, ~~~ n \in {\mathbb N},
$$
then, 
$$
\chi(S_{r_n}^{n-1}) \ge \frac{C_{m(n)}^{m(n)/2}}{C_{m(n)}^{p(n)}}.
$$}

\vskip+0.2cm

Translating Theorem 4 into a form of Theorem 3, we get 

\vskip+0.2cm

\noindent {\bf Theorem 5.} {\it Consider a sequence $ \{r_n\}_{n=1}^{\infty} $, where $ r_n > \frac{1}{2} $ for each $ n \in {\mathbb N} $. Let 
$ p(n) $ be the same as in Theorem 4. If 
$$
\frac{m(n)}{4} < p(n) < \frac{m(n)}{2} - \sqrt{\frac{m(n) \ln m(n)}{\kappa}}, ~~~ \kappa < 2, ~~~ n \in {\mathbb N},
$$
then, 
$$
\chi(S_{r_n}^{n-1}) > n+1, ~~~ \forall ~ n \ge n_0.
$$}

\vskip+0.2cm

The quality of Theorem 5 depends on the estimates for the function $ f(x) $. Determining the exact asymptotic behaviour of $ f(x) $ is a 
very hard problem of analytical number theory (see \cite{Erd}). 
As far as we know, the best upper estimate is $ f(x) = O\left(x^{0.525-\varepsilon}\right) $
with a so small $ \varepsilon > 0 $ that the authors did not care of it (see \cite{BHP}). However, it is conjectured that $ f(x) = O(\ln x) $
(see \cite{Cr}). The tightest lower bound is given in \cite{Sch} and \cite{R}, but it is sublogarithmic and apparently far enough from the 
truth. Using this information, we may derive

\vskip+0.2cm

\noindent {\bf Theorem 6.} {\it Assume that $ c_0 > 0 $ is such that $ f(x) \le c_0 x^{0.525} $ for every $ x $. Then, there exists a constant 
$ c_0' > 0 $ such that for any sequence of radii $ r_n $ satisfying the inequality
$$
r_n \ge \frac{1}{2} + \frac{c_0'}{n^{0.475}},
$$
we have the bound
$$
\chi(S_{r_n}^{n-1}) > n+1, ~~~ \forall ~ n \ge n_0.
$$}

\vskip+0.2cm

\noindent {\bf Theorem 7.} {\it Assume that $ c_1 > 0 $ is such that $ f(x) \le c_1 \ln x $ for every $ x $. Then, there exists a constant 
$ c_1' > 0 $ such that for any sequence of radii $ r_n $ satisfying the inequality
$$
r_n \ge \frac{1}{2} + c_1' \sqrt{\frac{\ln n}{n}},
$$
we have the bound
$$
\chi(S_{r_n}^{n-1}) > n+1, ~~~ \forall ~ n \ge n_0.
$$}

\vskip+0.2cm

So $ r_n > \frac{1}{2} $ may be quite close to the value $ \frac{1}{2} $, and, nevertheless, the chromatic numbers will exceed 
the Lov\'asz upper estimate. Finally, it is of interest for which sequences of $ r_n $, we do really have the bound    
$ \chi(S_{r_n}^{n-1}) \le n+1 $.
  
\vskip+0.2cm

\noindent {\bf Theorem 8.} {\it There exists a constant 
$ c_2 > 0 $ such that for any sequence of radii $ r_n $ satisfying the inequality
$$
r_n \le \frac{1}{2} +  \frac{c_2}{n},
$$
we have the bound
$$
\chi(S_{r_n}^{n-1}) \le n+1, ~~~ \forall ~ n \ge n_0.
$$}

\vskip+0.2cm

Further structure of the paper is as follows. In Section 3, we shall give proofs for Theorems 1 -- 4. Section 4 will be devoted to 
proving Theorems 5 -- 7. In Section 5, we shall discuss Theorem 8. In Section 6, some more comments and suggestions will be given. 
In particular, we shall exhibit more general {\it upper} estimates for $ \chi(S_{r_n}^{n-1}) $ than those in Theorem 8. 
In Section 7, we shall present a general scheme for obtaining better (and, in some sense, optimal) constants $ \gamma $ than those 
appearing in Theorems 1 and 2.

\section{Proofs of Theorems 1 -- 4}

Among Theorems 1 -- 3, Theorem 2 covers both Theorem 1 and Theorem 3. So we start by proving Theorem 2. 

Fix an $ r \in \left(\frac{1}{2}, \frac{1}{\sqrt{2}}\right) $ and an $ n \in {\mathbb N} $. Let $ m < n $ be the maximum natural 
number which is divisible by 4. Let us find $ a' $ from the relation 
$$
\frac{\sqrt{m}}{\sqrt{2m-2a'}} = r, ~~~ {\rm i.e.}, ~~~ a' = \frac{m (2r^2-1)}{2r^2}.
$$
Let $ p $ be the smallest prime number satisfying the inequality 
$$
p > \frac{m-a'}{4} = \frac{m}{8r^2}.
$$
Set 
$$
a = m - 4p < a'.
$$
Consider the following graph $ G=(V,E) $:
$$
V = \left\{{\bf x}=(x_1, \dots, x_{m}):~ x_i \in \left\{-\frac{1}{\sqrt{2m-2a}}, \frac{1}{\sqrt{2m-2a}}\right\}, ~
x_1 + \ldots + x_{m} = 0\right\},
$$
$$
E = \{\{{\bf x},{\bf y}\}: ~ {\bf x}, {\bf y} \in V, ~ l_2({\bf x},{\bf y})=1\}.
$$

Obviously $ V \subset S_{r'}^{m-1} $, where 
$$
r' = \frac{\sqrt{m}}{\sqrt{2m-2a}} < \frac{\sqrt{m}}{\sqrt{2m-2a'}} = r.
$$
If we use the standard notation $ \chi(G) $ for the chromatic number of $ G $ and $ \alpha(G) $ for its independence number, then 
we get 
$$ 
\chi(S_r^{n-1}) \ge \chi(S_{r'}^{m-1}) \ge \chi(G) \ge \frac{|V|}{\alpha(G)} = \frac{C_{m}^{m/2}}{\alpha(G)}.
$$

So we are led to estimate $ \alpha(G) $ from above. It is convenient to transform $ G = (V,E) $ into an $ H = (W,F) $:
$$                                                                                                       
W = \left\{{\bf x} \cdot \sqrt{2m-2a}:~ {\bf x} \in V\right\}, ~~~~
F = \left\{\{{\bf x},{\bf y}\}: ~ {\bf x}, {\bf y} \in W, ~ l_2({\bf x},{\bf y})=\sqrt{2m-2a}\right\}.
$$
Let us denote by $ ({\bf x},{\bf y}) $ the Euclidean scalar product of $ {\bf x} $ and $ {\bf y} $. Since for any $ {\bf x} \in W $, 
$ ({\bf x},{\bf x}) = m $, we may rewrite $ F $ as follows:
$$
F = \{\{{\bf x},{\bf y}\}: ~ {\bf x}, {\bf y} \in W, ~ ({\bf x},{\bf y})=a\}.
$$

Notice that for $ {\bf x}, {\bf y} \in W $, the quantity $ ({\bf x},{\bf y}) $ lies in the interval $ [-m,m] $ and is congruent to zero 
modulo 4. The last observation is due to the fact that $ m \equiv 0 \pmod 4 $ and every vector $ {\bf x} \in W $ contains an even number 
of negative coordinates. Also, 
$$
m - 8p < m - \frac{m}{r^2} < - m.
$$
Thus, for every $ {\bf x}, {\bf y} \in W $, 
$$
({\bf x},{\bf y}) \equiv m \pmod p  \Longleftrightarrow ({\bf x},{\bf y}) = m ~ {\rm or} ~ ({\bf x},{\bf y}) = a.
\eqno{(1)}
$$ 

Now, we are about to prove that $ \alpha(G) = \alpha(H) \le C_m^p $. Take an arbitrary 
$$ 
Q = \{{\bf x}_1, \dots, {\bf x}_s\} \subset W, ~~~ \forall~ i ~ \forall ~ j, ~~ ({\bf x}_i,{\bf x}_j) \neq a.
\eqno{(2)}
$$
In other words, $ Q $ is an independent set in $ H $. We have to show that $ s \le C_m^p $. For this purpose, we use the 
linear algebra method (see \cite{FW}, \cite{ABS}, \cite{BF}, \cite{Rai8}). 

To each vector $ {\bf x} \in W $ we assign a polynomial $ P_{{\bf x}} \in {\mathbb Z}/p{\mathbb Z}[y_1, \dots, y_m] $. First, we 
take 
$$
P_{{\bf x}}'({\bf y}) = \prod_{i \in I} (i - ({\bf x},{\bf y})), 
$$
where 
$$
I = \{0,1, \dots, p-1\} \setminus \{m \pmod p\}, ~~~ {\bf y} = (y_1, \dots, y_m),
$$
and so $ P_{{\bf x}}' \in {\mathbb Z}/p{\mathbb Z}[y_1, \dots, y_m] $.  Obviously, 
$$
\forall ~ {\bf x}, {\bf y} \in W ~~~~ P_{{\bf x}}'({\bf y}) \equiv 0 \pmod p \Longleftrightarrow ({\bf x},{\bf y}) \not \equiv m \pmod p.
\eqno{(3)}
$$

Second, we represent $ P_{{\bf x}}' $ as a sum of monomials. If a monomial has the form 
$$ 
y_{i_1}^{\alpha_{i_1}} \cdot \ldots \cdot y_{i_q}^{\alpha_{i_q}}, ~~~ \alpha_{i_1} > 0, \dots, \alpha_{i_q} > 0,
$$
then we replace it by 
$$ 
y_{i_1}^{\beta_{i_1}} \cdot \ldots \cdot y_{i_q}^{\beta_{i_q}},
$$
where $ \beta_{\nu} = 1 $, provided $ \alpha_{\nu} $ is odd, and $ \beta_{\nu} = 0 $, provided $ \alpha_{\nu} $ is even. 
Eventually, we get a polynomial $ P_{{\bf x}} $. It is worth noting that this polynomial does also satisfy property (3). 

It follows from properties (1), (2), and (3) that the polynomials 
$$
P_{{\bf x}_1}, \dots, P_{{\bf x}_s}
$$
assigned to the vectors of the set $ Q $ are linearly independent over $ {\mathbb Z}/p {\mathbb Z} $. It is also easy to see that 
the dimension of the space generated by    
$$
P_{{\bf x}_1}, \dots, P_{{\bf x}_s}
$$
does not exceed $ C_m^p $. Thus, $ s = |Q| \le C_m^p $ and, therefore, 
$$ 
\chi(S_r^{n-1}) \ge \frac{C_m^{m/2}}{C_m^p}. 
$$

Standard analytical tools (like Stirling's formula) together with $ p \sim \frac{m}{8r^2} $ give us, finally, the expected bound
$$ 
\chi(S_r^{n-1}) \ge \left(2 \left(\frac{1}{8r^2}\right)^{\frac{1}{8r^2}} 
\left(1-\frac{1}{8r^2}\right)^{1-\frac{1}{8r^2}}+\delta(n)\right)^n,
$$ 
which completes the proof of Theorems 1 -- 3. 

\vskip+0.2cm

The proof of Theorem 4 is now clear. We just reproduce the above argument with $ r_n $ instead of $ r $. The only thing one has 
to explain here is why we impose additional conditions on the value of a prime. Indeed, the inequality $ p(n) > \frac{m(n)}{4} $ 
is quite important, since property (1) becomes false without it. As for the inequality $ p(n) \le \frac{m(n)}{2} $, it is 
necessary to correctly estimate the independence number of our graph $ G $ by the quantity $ C_m^p $. Moreover, $ \chi(G) = 1 $, 
provided $ p(n) > \frac{m(n)}{2} $, and the result is trivial. Theorem 4 is proved. 

\section{Proofs of Theorems 5 -- 7}

\subsection{Proof of Theorem 5}

Set $ m = m(n) $, $ p = p(n) $. Since the function $ \frac{C_m^{m/2}}{C_m^p} $ is decreasing in $ p $, we just have to show that for 
$$
p = \left[\frac{m}{2}-\sqrt{\frac{m\ln m}{\kappa}}\right], 
$$
the inequality $ \frac{C_m^{m/2}}{C_m^p} > n+1 $ is true for large values of $ n $. We have 
$$
\frac{C_m^{m/2}}{C_m^p} = 
\frac{\left(\frac{m}{2}+1\right) \cdot \left(\frac{m}{2}+2\right) \cdot \ldots \cdot \left(\frac{m}{2}+\left(\frac{m}{2}-p\right)\right)}{
\frac{m}{2} \cdot \left(\frac{m}{2}-1\right) \cdot \ldots \cdot \left(\frac{m}{2}-\left(\frac{m}{2}-p-1\right)\right)} = 
$$
$$
=\frac{\left(1+\frac{2}{m}\right) \cdot \left(1+\frac{4}{m}\right) \cdot \ldots \cdot \left(1+\frac{m-2p}{m}\right)}{
\left(1-\frac{2}{m}\right) \cdot \left(1-\frac{4}{m}\right) \cdot \ldots \cdot \left(1-\frac{m-2p-2}{m}\right)} \sim
e^{\frac{(m-2p)^2}{2m}} \ge e^{\frac{2\ln m}{\kappa}} = m^{\frac{2}{\kappa}}.
$$
By a condition of Theorem 5, $ \kappa < 2 $. Thus, 
$$
m^{\frac{2}{\kappa}}(1+o(1)) > n+1, ~~~ \forall~ n \ge n_0.
$$
Theorem 5 is proved.

\subsection{Proof of Theorem 6}

We just have to show that for our choice of $ r_n $, 
$$
p = \frac{m}{8r_n^2} + f\left(\frac{m}{8r_n^2}\right) < \frac{m}{2} - \sqrt{\frac{m\ln m}{\kappa}}, 
$$
provided $ \kappa < 2 $ is a constant and $ n $ is large enough. 

Indeed, assume that $ c_0' $ is large (say, $ c_0' > c_0 $). Then,
$$
p \le \frac{m}{8\left(\frac{1}{2}+\frac{c_0'}{n^{0.475}}\right)^2} + c_0 
\left(\frac{m}{8\left(\frac{1}{2}+\frac{c_0'}{n^{0.475}}\right)^2}\right)^{0.525} <
$$
$$
< \frac{m}{8\left(\frac{1}{4}+\frac{c_0'}{n^{0.475}}\right)} + c_0 
\left(\frac{m}{8\left(\frac{1}{4}+\frac{c_0'}{n^{0.475}}\right)}\right)^{0.525} =
$$ 
$$
= \frac{m}{2}\left(1-\frac{4c_0'}{n^{0.475}} + o\left(\frac{1}{\sqrt{n}}\right)\right) + c_0
\left(\frac{m}{2}\left(1-\frac{4c_0'}{n^{0.475}} + o\left(\frac{1}{\sqrt{n}}\right)\right)\right)^{0.525}.
$$ 
For any sufficiently large value of $ n $, the last quantity is bounded from above by 
$$
\frac{m}{2} - c_0' m^{0.525} + c_0 m^{0.525} = \frac{m}{2} - c_0'' m^{0.525}, ~~~ c_0'' > 0.
$$
Obviously, for any $ n \ge n_0 $, 
$$
\frac{m}{2} - c_0'' m^{0.525} < \frac{m}{2} - \sqrt{\frac{m\ln m}{\kappa}}.
$$
Theorem 6 is proved.

\subsection{Proof of Theorem 7}

Let us briefly write down a series of inequalities similar to those in 4.2:
$$
p \le \frac{m}{8\left(\frac{1}{2}+c_1'\sqrt{\frac{\ln n}{n}}\right)^2} + c_1 
\ln \left(\frac{m}{8\left(\frac{1}{2}+c_1'\sqrt{\frac{\ln n}{n}}\right)^2}\right) <
$$
$$
< \frac{m}{2}\left(1-4c_1'\sqrt{\frac{\ln n}{n}} + o\left(\frac{1}{n^{3/2}}\right)\right) + c_1 \ln \left(
\frac{m}{2}\left(1-4c_1'\sqrt{\frac{\ln n}{n}} + o\left(\frac{1}{n^{3/2}}\right)\right)\right) <
$$
$$
< \frac{m}{2} - c_1' \sqrt{m \ln m},
$$
and we are done. 
 
\section{Proof of Theorem 8}

Let us take $ S_{1/2}^{n-1} $ and divide it into $ n+1 $ parts of smallest possible diameters. To this end, we inscribe a 
regular $ n $-simplex $ \Delta^n $ into $ S_{1/2}^{n-1} $ and consider multidimensional polyhedral cones
$ C_1, \dots, C_{n+1} $ with common vertex 
at the center of $ S_{1/2}^{n-1} $ and coming through the $ (n-1) $-faces of $ \Delta^n $. Obviously, 
$$
S_{1/2}^{n-1} = (S_{1/2}^{n-1} \cap C_1) \cup \ldots \cup (S_{1/2}^{n-1} \cap C_{n+1}).
\eqno{(4)}
$$
In principle, it is a good exercise in multidimensional geometry to prove that for any $ i $, 
$$
{\rm diam}~(S_{1/2}^{n-1} \cap C_i) = 1 - \Theta\left(\frac{1}{n}\right).
$$
It follows immediately from this observation that we may inflate $ S_{1/2}^{n-1} $ at most
$$
\frac{1}{1 - \Theta\left(\frac{1}{n}\right)} = 1 + \Theta\left(\frac{1}{n}\right)
$$
times in order to get a partition of the resulting sphere into parts of diameter not exceeding 1. Thus, for a constant 
$ c_2 > 0 $, we have an appropriate coloring of $ S_{1/2+c_2/n}^{n-1} $, which completes the proof of Theorem 8. 

\vskip+0.2cm

Apparently, in \cite{Lov}, the same construction was proposed. However, the author assumed that the diameter of any part in the 
corresponding partition is attained on the sides of a regular $ n $-simplex $ \Delta^n $. This is true only for $ n = 2 $. 
Already in $ {\mathbb R}^3 $, the diameter of a part is $ \sqrt{\frac{3+\sqrt{3}}{6}} = 0.888... $, which is not the length of a side 
of a tetrahedron inscribed into $ S_{1/2}^2 $.

\section{Comments and upper bounds}

First of all, it is worth noting that there is still a certain gap between the estimates 
$$
r_n \ge \frac{1}{2} + c_1' \sqrt{\frac{\ln n}{n}}
\eqno{(5)}
$$
and 
$$
r_n \le \frac{1}{2} +  \frac{c_2}{n}. 
\eqno{(6)}
$$
Removing this gap could be a good problem. As for (6), it cannot be enlarged by any refinement of the techniques of the previous section. 
The point is that the partition (4) is best possible: for any other decomposition of $ S_{1/2}^{n-1} $ into $ n+1 $ parts, 
there exists a part whose diameter is not less than each of the diameters $ {\rm diam}~(S_{1/2}^{n-1} \cap C_i) $. Of course it is not 
necessary to divide a sphere into parts with diameters strictly smaller than 1; we just need to cut it in such a way that no part would 
contain a pair of points at the unit distance. However, we do not know such a partition. Perhaps it is easier to improve (5). One should 
combine linear algebra of Section 3 with some additional ideas.

\vskip+0.2cm

Let us say a few words about general upper estimates for $ \chi(S_{r_n}^{n-1}) $. The simplest observation here is that 
$$
\chi(S_{r_n}^{n-1}) \le \chi({\mathbb R}^n) \le (3+o(1))^n ~~~ ({\rm cf.}~ {\rm Introduction}).
\eqno{(7)}
$$
Thus, for constant values $ r > \frac{1}{2} $ (as in Theorems 1 -- 3), we already get the order of magnitude for any
quantity $ \log \chi(S_r^{n-1}) $. 

In \cite{Rog}, C.A. Rogers proved that any sphere of radius $ r $ in $ {\mathbb R}^n $ can be covered by 
$ \left(\frac{r}{\rho} + o(1)\right)^n $ spheres of radius $ \rho < r $. In our case, this means that 
$$
\chi(S_{r_n}^{n-1}) \le (2r_n+o(1))^n.
$$
If $ r_n < 3/2 $, then this bound is better than that in (7). 

More precisely, Rogers' estimate is as follows: {\it there is an absolute constant $ c > 0 $ such that, if $ r > \frac{1}{2} $ and
$ n \ge 9 $, any $ n $-dimensional spheres of radius $ r $ can be covered by less than $ cn^{5/2} (2r)^n $ spheres of radius $ \frac{1}{2} $.} 
A so precise formulation is unuseful when $ r $ is a constant, but coming again to $ r_n \to \frac{1}{2} $ we may carefully apply this 
statement in order to obtain upper bounds like
$$ 
\chi(S_{r_n}^{n-1}) \le 2cn^{5/2} (2r_n)^n = \Theta\left(n^{5/2} (2r_n)^n\right).
\eqno{(8)}
$$
Here the factor 2 is due to the fact that $ \chi(S_{1/2}^{n-1}) = 2 $. One should not forget that if, for example,  
$ r_n = \frac{1}{2} + \Theta\left(\frac{1}{n}\right) $, then $ (2r_n)^n = \Theta(1) $, so that estimate (8) is very good. 

It is possible to evaluate even more sophisticated bounds for $ \chi(S_{r_n}^{n-1}) $, but this is not so interesting. 

\section{A possible way for improving Theorem 2}

\subsection{Statements of the results}

Fix again an $ r > \frac{1}{2} $. Let $ m = m(n) < n $ for every $ n $ and $ m \sim n $ for $ n \to \infty $. 
Assume that $ t = t(n) \in {\mathbb N} $, 
$$ 
b_1 = b_1(n) \in {\mathbb Z}, ~~\dots, ~~b_t = b_t(n) \in {\mathbb Z}, 
$$ 
$$ 
l_1 = l_1(n) \in {\mathbb N},~~ \dots,~~ l_t = l_t(n) \in {\mathbb N}, ~~~ l_1 + \ldots + l_t = m.
$$ 
Consider
$$
V = V(n) = \{{\bf x} = (x_1, \dots, x_{m}):~ x_i \in \{b_1, \dots, b_t\}, ~ |\{i:~ x_i = b_j\}| = l_j, ~ j = 1, \dots, t\}.
$$
Let $ d = d(n) $ be the maximum natural number such that for any 
$ {\bf x}, {\bf y} \in V $, we have $ ({\bf x},{\bf y}) \equiv 0 \pmod d $. 
Note that $ V $ is an obvious analog of the set $ W $ from Section 3, where $ d $ was equal to 4. 
Set 
$$
\overline{s} = \overline{s}(n) = \max_{{\bf x},{\bf y} \in V} ({\bf x},{\bf y}), ~~
\underline{s} = \underline{s}(n) = \min_{{\bf x},{\bf y} \in V} ({\bf x},{\bf y}).
$$
Find $ a' = a'(n) $ from the relation
$$
\frac{\sqrt{\overline{s}}}{\sqrt{2\overline{s}-2a'}} = r.
$$
Define $ p = p(n) $ as the minimum prime number satisfying the inequality
$$
p > \frac{\overline{s}-a'}{d}.
$$
Finally, we choose $ a = a(n) $ from the condition
$$
p = \frac{\overline{s}-a}{d}, ~~~ {\rm i.e.,} ~~~ a = \overline{s}-dp < a'.
$$
We get the following theorem.

\vskip+0.2cm

\noindent {\bf Theorem 9.} {\it If $ a > \underline{s} $ and $ \overline{s}-2dp < \underline{s} $, then
$$
\chi(S_r^{n-1}) \ge \frac{L}{M},
$$
where 
$$
L = \frac{m!}{l_1! \cdot \ldots \cdot l_t!}, ~~~~ M = \sum_{(s_1, \dots, s_t) \in {\cal A}} \frac{m!}{s_1! \cdot \ldots \cdot s_t!},
$$
$$
{\cal A} = \{(s_1, \dots, s_t):~ s_i \in {\mathbb N} \cup \{0\}, ~ s_1 + \ldots + s_t = m, ~ s_1 + 2s_2 + \ldots + (t-1)s_{t-1} \le p-1\}.
$$}

\vskip+0.2cm

In Theorem 9 we optimize over the parameters
$ t $, $ b_1, \dots, b_t $, and $ l_1, \dots, l_t $. This optimization can be a bit simpler, provided we suppose that 
$ l_i \sim l_i^0 n $, where $ l_i^0 \in (0,1) $. Actually this does not substantially change results. In our case, we get

\vskip+0.2cm

\noindent {\bf Corollary.} {\it The estimate holds
$$
\chi(S_r^{n-1}) \ge \left(\frac{L_0}{M_0}+o(1)\right)^n,
$$
where 
$$
L_0 = e^{-l_1^0 \ln l_1^0 - \ldots - l_t^0 \ln l_t^0}, ~~~~ M_0 = \max_{(s_1^0, \dots, s_t^0) \in {\cal A}_0} 
e^{-s_1^0 \ln s_1^0 - \ldots - s_t^0 \ln s_t^0},
$$
$$
{\cal A}_0 = \left\{\left(s_1^0, \dots, s_t^0\right):~ s_i^0 \in (0,1), ~ s_1^0 + \ldots + s_t^0 = 1, ~ s_1^0 + 2s_2^0 + \ldots + (t-1)s_{t-1}^0 \le 
\frac{p}{n}\right\}.
$$}

\vskip+0.2cm

We shall prove Theorem 9 in \S 7.2. Corollary can be easily derived from Theorem 9 using Stirling's formula and other standard tools of analysis. 

In this paper, we shall not evaluate optimization from Corollary. Here we only cite the papers \cite{Rai3}, \cite{Rai7}, in which similar 
optimization procedures were carefully realized.

\subsection{Proof of Theorem 9}

Let us start by noting that all the parameters in Theorem 9 are chosen to generalize the approach that we used in Section 3. We have already 
mentioned that the quantity $ d $ plays the role of the number 4 in the corresponding argument. Almost all the other notations are also completely 
parallel to those appearing in Section 3. Here only $ m $ should be replaced by $ \overline{s} $, and we just consider $ V $ as an analog to $ W $, 
without introducing two similar sets $ V $ and $ W $ as it was done in Section 3. 

Set $ G = (V,E) $ with
$$
E = \{\{{\bf x},{\bf y}\}: ~ {\bf x}, {\bf y} \in V, ~ ({\bf x},{\bf y})=a\}.
$$
We think it is now obvious that 
$$
\chi(S_r^{n-1}) \ge \chi(G) \ge \frac{|V|}{\alpha(G)} = \frac{L}{\alpha(G)}.
$$
So it remains to prove that $ \alpha(G) \le M $. This is done by the same linear algebra method as in Section 3. 

Indeed, by the conditions of Theorem 9, we have, for every $ {\bf x}, {\bf y} \in V $, 
$$
({\bf x},{\bf y}) \equiv \overline{s} \pmod p  \Longleftrightarrow ({\bf x},{\bf y}) = \overline{s} ~ {\rm or} ~ ({\bf x},{\bf y}) = a.
\eqno{(1')}
$$ 
Take an arbitrary 
$$ 
Q = \{{\bf x}_1, \dots, {\bf x}_s\} \subset V, ~~~ \forall~ i ~ \forall ~ j, ~~ ({\bf x}_i,{\bf x}_j) \neq a.
\eqno{(2')}
$$
We are about to show that $ s \le M $.

To each vector $ {\bf x} \in V $ we assign a polynomial $ P_{{\bf x}} \in {\mathbb Z}/p{\mathbb Z}[y_1, \dots, y_m] $. First, we 
take 
$$
P_{{\bf x}}'({\bf y}) = \prod_{i \in I} (i - ({\bf x},{\bf y})), 
$$
where 
$$
I = \{0,1, \dots, p-1\} \setminus \{\overline{s} \pmod p\}, ~~~ {\bf y} = (y_1, \dots, y_m),
$$
and so $ P_{{\bf x}}' \in {\mathbb Z}/p{\mathbb Z}[y_1, \dots, y_m] $.  Obviously, 
$$
\forall ~ {\bf x}, {\bf y} \in W ~~~~ P_{{\bf x}}'({\bf y}) \equiv 0 \pmod p \Longleftrightarrow ({\bf x},{\bf y}) \not \equiv \overline{s} \pmod p.
\eqno{(3')}
$$

Second, we represent $ P_{{\bf x}}' $ as a sum of monomials. We use the fact that 
$$
(y_i-b_1)\cdot(y_i-b_2)\cdot \ldots \cdot (y_i-b_t)=0,
$$
for any $ {\bf y} \in V $. So we get a polynomial $ P_{{\bf x}} $ of degree $ < t $. 
It is worth noting that this polynomial does also satisfy property $ (3') $. 

It follows from properties $ (1'), (2'), $ and $ (3') $ that the polynomials 
$$
P_{{\bf x}_1}, \dots, P_{{\bf x}_s}
$$
assigned to the vectors of the set $ Q $ are linearly independent over $ {\mathbb Z}/p {\mathbb Z} $. Now it is easy to see that 
the dimension of the space generated by    
$$
P_{{\bf x}_1}, \dots, P_{{\bf x}_s}
$$
does not exceed $ M $. Thus, $ s = |Q| \le M $ and, therefore, Theorem 9 is proved.

\end{large}

\end{document}